\documentclass{article}
\usepackage{amsthm}
\usepackage{amssymb}
\usepackage{amsxtra}

\newcommand{\Dst}[1]{{\displaystyle #1}}
\newtheorem{thm}{Theorem}[section]
\newtheorem{lemma}{Lemma}[section]

\newtheorem{conj}{Conjecture}[section]

\title{Pair Correlation of the zeros of the Riemann zeta function
in longer ranges}
\author{Tsz Ho Chan}

\begin{document}
\maketitle
\begin{abstract}
In this paper, we study a more general pair correlation function,
$F_h(x,T)$, of the zeros of the Riemann zeta function. It provides
information on the distribution of larger differences between the
zeros.
\end{abstract}
\section{Introduction}

First of all, we assume the Riemann Hypothesis on the Riemann zeta
function $\zeta(s)$ throughout this paper; $\rho = \frac{1}{2} +
i\gamma$ denotes a non-trivial zero of the Riemann zeta function.

In the early 1970s, Hugh Montgomery considered the pair
correlation function
\begin{equation*}
F(x,T) = \sum_{0 < \gamma, \gamma' \leq T}
         x^{i(\gamma - {\gamma}')} w(\gamma - {\gamma}')
\mbox{ with } w(u) = \frac{4}{4 + u^2}.
\end{equation*}
Here the sum is a double sum over the imaginary parts of the
non-trivial zeros of $\zeta(s)$. He proved in [\ref{M}] that, as
$T \rightarrow \infty$,
$$F(x,T) \sim \frac{T}{2 \pi} \log{x} + \frac{T}{2 \pi x^2} \log^2{T}$$
for $1 \leq x \leq T$ (actually he only proved for $1 \leq x \leq
o(T)$ and the full range was done by Goldston [\ref{G}]). He
conjectured that
$$F(x,T) \sim \frac{T}{2 \pi} \log{T}$$
for $T \leq x$ which is known as the Strong Pair Correlation
Conjecture. From this, one has the (Weak) Pair Correlation
Conjecture:
$$\mathop{\sum_{0 < \gamma, \gamma' \leq T}}_{0 < \gamma - \gamma' \leq 2 \pi
\alpha / \log{T}} 1 \sim \frac{T}{2 \pi}\log{T} \int_{0}^{\alpha}
1 - \Bigl(\frac{\sin{\pi u}}{\pi u}\Bigr)^2 du.$$ which draws
connections with random matrix theory.

The author studied these further in his thesis [\ref{C1}] (see
also [\ref{C2}] and [\ref{C3}]) and derived more precise
asymptotic formulas for $F(x,T)$ when $x$ is in various ranges
under the Twin Prime Conjecture {\bf TPC} (see section 4). In the
present paper, we generalize $F(x,T)$ further to
$$F_h(x,T) = \sum_{0 < \gamma, \gamma' \leq T} \cos{\bigl((\gamma -
\gamma' - h) \log{x}\bigr)} w(\gamma - \gamma' - h).$$ Note that
$F_h(x,T) = F_{-h}(x,T)$ and $F_0(x,T) = F(x,T)$. This leads to a
better understanding of the distribution of larger differences
between the zeros. Our main results are the following theorems:
Here and throughout the paper, $\tilde{h} = |h|+1$.

\begin{thm}
\label{theorem1} For $1 \leq x \leq \frac{T}{\log{T}}$,
\begin{equation*}
\begin{split}
F_h(x, T) =& \frac{T}{2 \pi} \Bigl[\frac{4 \cos{(h \log{x})}}{4 +
h^2} \log{x} - \frac{8h \sin{(h \log{x})}}{(4 + h^2)^2} \Bigr] \\
& + \frac{T}{2 \pi x^2} \Bigl[ \Bigl(\log{\frac{T}{2 \pi}}\Bigr)^2
- 2\log{\frac{T}{2 \pi}} \Bigr] + O(x \log{x}) + O\Bigl(
\frac{\tilde{h} T}{x^{1/2 - \epsilon}} \Bigr).
\end{split}
\end{equation*}
\end{thm}

\begin{thm}
\label{theorem2} Assume {\bf TPC}. For $M \geq 3$ and
$\frac{T}{\log^M{T}} \leq x$,
\begin{equation*}
\begin{split}
& F_h(x,T) \\
=& \frac{T}{\pi} \Bigl[\frac{2 \cos{(h\log{x})}}{4 + h^2} \log{x}
- \frac{4h \sin{(h \log{x})}}{(4 + h^2)^2} \Bigr] \\
+& \frac{T}{\pi} \int_{1}^{\infty} \Bigl[ - \frac{2
\cos{(h\log{x})}} {4+h^2} \frac{1}{y} - \frac{4f(y)}{y^2}
\cos{(h\log{x})} + G_1(y) + G_2(y)\Bigr] \frac{\sin{\frac{Ty}{x}}}
{\frac{Ty}{x}} dy\\
-& \frac{x}{\pi} \int_{0}^{T/x} \frac{\sin{u}}{u} du
\Bigl[\frac{3\cos{(h\log{x})}}{9+h^2} +
\frac{h\sin{(h\log{x})}}{9+h^2}\Bigr] \\
-& \frac{x}{\pi} \int_{0}^{T/x} \frac{\sin{u}}{u} du
\Bigl[\frac{\cos{(h\log{x})}}{1+h^2} -
\frac{h\sin{(h\log{x})}}{1+h^2}\Bigr] \\
+& \frac{T}{\pi} \sum_{k=1}^{\infty} \frac{{\mathfrak S}(k)}{k^2}
\int_{0}^{1} y \cos{(h\log{\frac{kx}{y}})}
\frac{\sin{\frac{Ty}{x}}} {\frac{Ty}{x}} dy \\
+& O\Bigl(\tilde{h}\frac{x^{1+6\epsilon}}{T}\Bigr) + O(\tilde{h}
x^{1/2 + 7\epsilon}) + O\Bigl(\tilde{h}
\frac{x^2}{T^{2-2\epsilon}}\Bigr) +
O\Bigl(\frac{\tilde{h}T}{\log^{M-2}{T}}\Bigr).
\end{split}
\end{equation*}
where $G_1(y)$ and $G_2(y)$ are defined in Lemma \ref{lemma4.2}.
\end{thm}

\begin{thm}
\label{theorem3} Assume {\bf TPC}. For $M \geq 3$ and
$\frac{T}{\log^M{T}} \leq x \leq T$,
$$F_h(x,T) = \frac{T}{\pi} \Bigl[\frac{2 \cos{(h\log{x})}}{4+h^2}
\log{x} - \frac{4h \sin{(h\log{x})}}{(4+h^2)^2} \Bigr] +
O(\tilde{h} x) + O\Bigl(\frac{\tilde{h}T}{\log^{M-2}{T}}\Bigr).$$
\end{thm}

\begin{thm}
\label{theorem4} Assume {\bf TPC}. For $M \geq 3$ and $T \leq x
\leq T^{2-29\epsilon}$,
$$F_h(x,T) = \frac{T}{2\pi} \log{\frac{T}{2\pi e}} \Bigl[\frac{4
\cos{(h\log{x})}}{4 + h^2}\Bigr] + O\Bigl(\tilde{h} T
\bigl(\frac{T}{x}\bigr)^{1/2-\epsilon}\Bigr) +
O\Bigl(\frac{\tilde{h} T} {\log^{M-2}{T}}\Bigr).$$
\end{thm}

For real $\alpha$, let $F_h(\alpha) := (\frac{T}{2\pi}
\log{T})^{-1} F_h(T^{\alpha}, T)$. Then $F_h(\alpha) =
F_h(-\alpha)$. Based on the above theorems, one may make the
following
\begin{conj}
\label{conj1} For any arbitrary large $A$, as $T \rightarrow
\infty$,
\[F_h(\alpha) = \left\{ \begin{array}{ll}
(1 + o(1))T^{-2\alpha} \log{T} + \alpha \frac{4\cos{(h\log{T}
\alpha)}}{4+h^2} + o(1), & \mbox{ if }\, 0 \leq \alpha \leq 1, \\
\frac{4\cos{(h\log{T} \alpha)}}{4+h^2} + o(1), & \mbox{ if }\, 1
\leq \alpha \leq A. \end{array} \right.
\]
\end{conj}
By convolving $F_h(\alpha)$ with an appropriate kernel
$\hat{r}(\alpha)$,
\begin{equation}
\label{0} \Bigl(\frac{T}{2\pi} \log{T}\Bigr)^{-1} \sum_{0 <
\gamma, \gamma' \leq T} r\Bigl((\gamma - \gamma' - h)
\frac{\log{T}}{2\pi} \Bigr) w(\gamma -\gamma' -h) =
\int_{-\infty}^{+\infty} F_h(\alpha) \hat{r}(\alpha)\, d\alpha
\end{equation}
where $\hat{r}(\alpha) = \int_{-\infty}^{\infty} r(u) e^{-2\pi i
\alpha u} du$ for even $r(u)$ only. Conjecture \ref{conj1} and
(\ref{0}) leads to
\begin{conj}
\label{conj2} For fixed $\alpha > 0$,
$$\Bigl(\frac{T}{2\pi}\log{T}\Bigr)^{-1} \mathop{\sum_{0 < \gamma
\neq \gamma' \leq T}}_{|\gamma - \gamma' - h| \leq 2\pi\alpha /
\log{T}} 1 \sim \int_{-\alpha + h\log{T}/(2 \pi)}^{\alpha +
h\log{T}/(2 \pi)} 1 - \frac{4}{4+h^2} \Bigl(\frac{\sin{\pi u}}{\pi
u}\Bigr)^2 du.$$
\end{conj}
\begin{conj}
For $0 < \alpha < \beta \ll \log{T}$,
$$\Bigl(\frac{T}{2\pi}\log{T}\Bigr)^{-1} \mathop{\sum_{0 < \gamma
\neq \gamma' \leq T}}_{2\pi\alpha/\log{T} \leq \gamma - \gamma'
\leq 2\pi\beta / \log{T}} 1 \sim \int_{\alpha}^{\beta} 1 -
\frac{1}{1 + (\pi u/\log{T})^2} \Bigl(\frac{\sin{\pi u}}{\pi
u}\Bigr)^2 du.$$
\end{conj}
\section{Some Lemmas}
\begin{lemma}
\label{lemma1.1}
\begin{equation}
\label{1.1}
\begin{split}
2 \sum_{\gamma} \frac{x^{i (\gamma - t)}}{1 + (t - \gamma)^2} =&
-\frac{1}{x} \sum_{n \leq x} \frac{\Lambda(n)}{n^{-1/2+i t}} - x
\sum_{n > x} \frac{\Lambda(n)}{n^{3/2 + it}} + \frac{x^{1/2 -
it}}{\frac{1}{2} + i t} + \frac{x^{1/2 - i t}}{\frac{3}{2} - i t}
\\
&+ \frac{\log{\tau}}{x} + \frac{1}{x} \Bigl[
\frac{\zeta'}{\zeta}\Bigl(\frac{3}{2} - i t \Bigr) - \log{2
\pi}\Bigr] + O\Bigl(\frac{1}{x \tau}\Bigr)
\end{split}
\end{equation}
where the sum is over all the imaginary parts of the zeros of the
Riemann zeta function, and $\tau = |t| + 2$. $\Lambda(n)$ is von
Mangoldt's lambda function.
\end{lemma}

Proof: This is Lemma 2.2 in [\ref{C2}].

Write (\ref{1.1}) as $\mbox{Left}(x,t) = \mbox{Right}(x,t)$. Let
\begin{equation*}
\begin{split}
P(x, T) =& \frac{1}{x} \sum_{n \leq x} \frac{\Lambda(n)}{n^{-1/2+i
t}} + x \sum_{n > x} \frac{\Lambda(n)}{n^{3/2 + it}} -
\frac{x^{1/2 - it}}{\frac{1}{2} + i t} - \frac{x^{1/2 - i
t}}{\frac{3}{2} - i t}, \\
Q(x, T) =& \frac{\log{\tau}}{x}, \\
R(x, T) =& \frac{1}{x} \Bigl[
\frac{\zeta'}{\zeta}\Bigl(\frac{3}{2} - i t \Bigr) - \log{2
\pi}\Bigr], \\
S(x, T) =& O\Bigl(\frac{1}{x \tau}\Bigr).
\end{split}
\end{equation*}

\begin{lemma}
\label{lemma1.2}
\begin{equation*}
\begin{split}
&\int_{0}^{T} |\mbox{Left}(x,t) + \mbox{Left}(x,t-h)|^2 dt \\
=& 2\pi F(x,T) + 2\pi F(x,T-h) + 4\pi F_h(x,T) + O(\log^3{T}) +
O(h \log^2{h}).
\end{split}
\end{equation*}
\end{lemma}

Proof: This follows from page $188$ of Montgomery [\ref{M}] and
the fact that $F(x,T) \ll T \log^2{T}$.

\begin{lemma}
\label{lemma1.3}
For $x \geq 1$,
\begin{equation*}
\begin{split}
&\int_{0}^{T} |\mbox{Right}(x,t) + \mbox{Right}(x,t-h)|^2 dt \\
=&  \int_{0}^{T} |P(x,t) + P(x,t-h)|^2 dt + \frac{4 T}{x^2} \Bigl[
\Bigl(\log{\frac{T}{2 \pi}}\Bigr)^2 - 2\log{\frac{T}{2 \pi}} \\
&+ \Bigl(\frac{1}{2} \sum_{n = 1}^{\infty} \frac{\Lambda^2(n) (1 +
\cos{(h\log{n})})}{n^3} + 2\Bigr) \Bigr] + O(\tilde{h} \log^2{T}).
\end{split}
\end{equation*}
\end{lemma}

Proof: This is similar to the proof of Theorem 3.1 in [\ref{C2}].

\begin{lemma}
\label{lemma1.4}
For $x \geq 1$,
\begin{equation*}
\begin{split}
& 4 \pi F_h(x,T) \\
=& \int_{0}^{T} |P(x,t) + P(x,t-h)|^2 dt - \int_{0}^{T}
|P(x,t)|^2 dt - \int_{0}^{T} |P(x,t-h)|^2 dt \\
& + \frac{2 T}{x^2} \Bigl[ \Bigl(\log{\frac{T}{2 \pi}}\Bigr)^2 -
2\log{\frac{T}{2 \pi}} + \Bigl(\sum_{n = 1}^{\infty}
\frac{\Lambda^2(n) \cos{(h\log{n})}}{n^3} + 2\Bigr) \Bigr] +
O(\tilde{h} \log^3{T}).
\end{split}
\end{equation*}
\end{lemma}

Proof: It follows from Lemma \ref{lemma1.2} and Lemma
\ref{lemma1.3} as well as their special cases when $h = 0$.

\begin{lemma}
\label{lemma1.5}
For any sequence of complex numbers
$\Dst{\{a_n\}_{n=1}^{\infty}}$ with $\Dst{\sum_{n=1}^{\infty} n
|a_n|^2 < \infty}$,
$$\int_{0}^{T} \Bigl|\sum_{n=1}^{\infty} a_n n^{-it} \Bigr|^2 dt =
\sum_{n=1}^{\infty} |a_n|^2 \bigl(T + O(n)\bigr).$$
\end{lemma}

Proof: This is Parseval's identity for Dirichlet series. See
[\ref{MV}].

\begin{lemma}
\label{lemma1.6}
\begin{equation*}
\begin{split}
\sum_{n \leq x} \Lambda^2(n) n =& \frac{1}{2} x^2\log{x} -
\frac{1}{4} x^2 + O(x^{1/2+\epsilon}). \\
\sum_{n > x} \frac{\Lambda^2(n)}{n^3} =& \frac{1}{2}
\frac{\log{x}}{x^2} + \frac{1}{4 x^2} +
O\Bigl(\frac{1}{x^{5/2-\epsilon}}\Bigr).
\end{split}
\end{equation*}
\end{lemma}

Proof: Use partial summation and the prime number theorem.

\begin{lemma}
\label{lemma1.7} For any real $a$ and $b$ not both zero,
\begin{equation*}
\begin{split}
\int e^{a x} \sin{b x}\, dx =& \frac{a}{a^2+b^2} e^{a x} \sin{b
x} - \frac{b}{a^2+b^2} e^{a x} \cos{b x}. \\
\int e^{a x} \cos{b x}\, dx =& \frac{a}{a^2+b^2} e^{a x} \cos{b
x} + \frac{b}{a^2+b^2} e^{a x} \sin{b x}. \\
\int x e^{a x} \sin{b x}\, dx =& \Bigl[\frac{a x}{a^2+b^2} -
\frac{a^2-b^2}{(a^2+b^2)^2} \Bigr] e^{a x} \sin{b x} \\
&- \Bigl[\frac{b x}{a^2+b^2} - \frac{2a b}{(a^2+b^2)^2} \Bigr]
e^{a x} \cos{b x}. \\
\int x e^{a x} \cos{b x}\, dx =& \Bigl[\frac{a x}{a^2+b^2} -
\frac{a^2-b^2}{(a^2+b^2)^2} \Bigr] e^{a x} \cos{b x} \\
&+ \Bigl[\frac{b x}{a^2+b^2} - \frac{2a b}{(a^2+b^2)^2} \Bigr]
e^{a x} \sin{b x}.
\end{split}
\end{equation*}
\end{lemma}

Proof: One can use $\Dst{\int e^{(a+ib)x} dx}$, $\Dst{\int
e^{(a-ib)} dx}$, $\Dst{\int x e^{(a+ib)x} dx}$ and $\Dst{\int x
e^{(a-ib)x} dx}$ which are simple to compute.

\begin{lemma}
\label{lemma1.8}
\begin{equation*}
\begin{split}
\frac{1}{x^2} \sum_{n \leq x} \Lambda^2(n) n \cos{(h\log{n})} =&
\frac{2\cos{(h\log{x})}}{4+h^2} \log{x} + \frac{h^2-4}{(4+h^2)^2}
\cos{(h\log{x})} \\
&+ \frac{h\sin{(h\log{x})}}{4+h^2} \log{x} - \frac{4h}{(4+h^2)^2}
\sin{(h\log{x})} \\
&+ O\Bigl(\frac{\tilde{h}}{x^{1/2-\epsilon}}\Bigr). \\
x^2 \sum_{n > x} \frac{\Lambda^2(n)}{n} \cos{(h\log{n})} =&
\frac{2\cos{(h\log{x})}}{4+h^2} \log{x} - \frac{h^2-4}{(4+h^2)^2}
\cos{(h\log{x})} \\
&- \frac{h\sin{(h\log{x})}}{4+h^2} \log{x} - \frac{4h}{(4+h^2)^2}
\sin{(h\log{x})} \\
&+ O\Bigl(\frac{\tilde{h}}{x^{1/2-\epsilon}}\Bigr).
\end{split}
\end{equation*}
\end{lemma}

Proof: We shall prove the first one. The other one is very
similar. Let $\Dst{A(x) = \frac{1}{x^2} \sum_{n \leq x}
\Lambda^2(n) n}$. By integration by parts and Lemma
\ref{lemma1.6},
\begin{equation*}
\begin{split}
& \frac{1}{x^2} \sum_{n \leq x} \Lambda^2(n) n \cos{(h\log{n})} \\
=& \frac{A(x)}{x^2} \cos{(h\log{x})} + \frac{h}{x^2} \int_{1}^{x}
A(u) \frac{\sin{(h\log{u})}}{u}\, du \\
=& \Bigl[\frac{1}{2} \log{x} - \frac{1}{4}\Bigr] \cos{(h\log{x})}
+ \frac{h}{x^2} \int_{1}^{x} \Bigl[\frac{1}{2} \log{u} -
\frac{1}{4}\Bigr] u \sin{(h\log{u})}\, du \\
&+ O\Bigl(\frac{\tilde{h}}{x^{1/2-\epsilon}}\Bigr) \\
=& \Bigl[\frac{1}{2} \log{x} - \frac{1}{4}\Bigr] \cos{(h\log{x})}
+ \frac{h}{x^2} \Bigl[ \frac{1}{2} \int_{0}^{\log{x}} v e^{2v}
\sin{hv}\, dv \\
&- \frac{1}{4} \int_{0}^{\log{x}} e^{2v} \sin{hv}\, dv \Bigr]+
O\Bigl(\frac{\tilde{h}}{x^{1/2-\epsilon}}\Bigr)
\end{split}
\end{equation*}
which gives the desired result after applying Lemma \ref{lemma1.7}
with $a=2$ and $b=h$, and some algebra.
\section{Proof of Theorem \ref{theorem1}}

First, note that
$$P(x,t) = \frac{1}{x^{1/2}} \Bigl[\sum_{n \leq x} \Lambda(n)
\Bigl(\frac{x}{n}\Bigr)^{-1/2+it} + \sum_{n > x} \Lambda(n)
\Bigl(\frac{x}{n}\Bigr)^{3/2+it} \Bigr] +
O\Bigl(\frac{x^{1/2}}{\tau}\Bigr).$$ Thus,
\begin{equation*}
\begin{split}
P(x,t) + P(x, t-h) =& \frac{1}{x^{1/2}} \Bigl[\sum_{n \leq x}
\Lambda(n) (1 + n^{ih}) \Bigl(\frac{x}{n}\Bigr)^{-1/2+it} \\
&+ \sum_{n > x} \Lambda(n) (1 + n^{ih})
\Bigl(\frac{x}{n}\Bigr)^{3/2+it} \Bigr] +
O\Bigl(\frac{x^{1/2}}{\tau}\Bigr)
\end{split}
\end{equation*}
So, the first integral in Lemma \ref{lemma1.4}
\begin{equation*}
\begin{split}
=& \frac{1}{x} \int_{0}^{T} \Bigl| \sum_{n \leq x} \Lambda(n) (1 +
n^{ih}) \Bigl(\frac{x}{n}\Bigr)^{-1/2+it} + \sum_{n > x}
\Lambda(n) (1 + n^{ih})
\Bigl(\frac{x}{n}\Bigr)^{3/2+it} \Bigr|^2 dt \\
&+ O\Bigl( \Bigl[ \sum_{n \leq x} \Lambda(n)
\Bigl(\frac{x}{n}\Bigr)^{-1/2} + \sum_{n > x} \Lambda(n)
\Bigl(\frac{x}{n}\Bigr)^{3/2}\Bigr] \int_{0}^{T} \frac{1}{\tau^2}
dt \Bigr) + O\Bigl(\int_{0}^{T} \frac{x}{\tau^4} dt \Bigr) \\
=& \frac{1}{x} \sum_{n \leq x} \Lambda^2(n) |1+n^{ih}|^2
\Bigl(\frac{x}{n}\Bigr)^{-1} (T + O(n)) \\
&+ \frac{1}{x} \sum_{n > x} \Lambda^2(n) |1+n^{ih}|^2
\Bigl(\frac{x}{n}\Bigr)^{3} (T + O(n)) + O(x) \\
=& \frac{2T}{x^2} \sum_{n \leq x} \Lambda^2(n) n \bigl(1 +
\cos{(h\log{n})} \bigr) + 2T x^2 \sum_{n > x}
\frac{\Lambda^2(n)}{n^3} \bigl(1 + \cos{(h\log{n})} \bigr) \\
&+ O(x \log{x})
\end{split}
\end{equation*}
Similarly (or by setting $h = 0$), each of the second and third
integral in Lemma \ref{lemma1.4}
\begin{equation*}
= \frac{T}{x^2} \sum_{n \leq x} \Lambda^2(n) n + T x^2 \sum_{n >
x} \frac{\Lambda^2(n)}{n^3} + O(x \log{x})
\end{equation*}
Therefore,
\begin{equation*}
\begin{split}
4 \pi F_h(x,T) =& 2T \Bigl[\frac{1}{x^2} \sum_{n \leq x}
\Lambda^2(n) n \cos{(h\log{n})} +  x^2 \sum_{n > x}
\frac{\Lambda^2(n)}{n^3} \cos{(h\log{n})} \Bigr] \\
&+ \frac{2 T}{x^2} \Bigl[ \Bigl(\log{\frac{T}{2 \pi}}\Bigr)^2 -
2\log{\frac{T}{2 \pi}}\Bigr] \\
&+ O\Bigl(\frac{T}{x^2}\Bigr) + O(\tilde{h} \log^3{T}) + O(x
\log{x}) \\
=& 2T \Bigl[\frac{4 \cos{(h\log{x})}}{4+h^2} \log{x} - \frac{8h
\sin{(h\log{x})}}{(4+h^2)^2} \Bigr] \\
&+ \frac{2 T}{x^2} \Bigl[ \Bigl(\log{\frac{T}{2 \pi}}\Bigr)^2 -
2\log{\frac{T}{2 \pi}}\Bigr] + O(x \log{x}) +
O\Bigl(\frac{\tilde{h}T}{x^{1/2-\epsilon}}\Bigr)
\end{split}
\end{equation*}
by Lemma \ref{lemma1.8}. The theorem follows after dividing
through by $4\pi$.
\section{Twin Prime Conjecture and smooth weight}

We shall use a quantitative form of the Twin Prime Conjecture {\bf
TPC} as follow: For any $\epsilon > 0$,
$$\sum_{n=1}^{N} \Lambda(n) \Lambda(n+d) = {\mathfrak S}(d) N + O(N^{1/2+
\epsilon}) \mbox{ uniformly in } |d| \leq N.$$ ${\mathfrak S}(d) =
2\prod_{p>2}\bigl(1-\frac{1}{(p-1)^2}\bigr) \prod_{p|d, p>2}
\frac{p-1}{p-2}$ if $d$ is even, and ${\mathfrak S}(d) = 0$ if $d$
is odd.

\bigskip

Let $K$ and $M$ be some large positive integers ($K$ may depend on
$\epsilon$). Set $U = \log^M{T}$ and $\Delta = 1/(2^K U)$. We
recall the smooth weight $\Psi_U(t)$ in [\ref{C3}] with:
\begin{enumerate}
\item support in $[-1/U, 1+1/U]$, \item $0 \leq \Psi_U (t) \leq
1$, \item $\Psi_U (t) = 1$ for $1/U \leq t \leq 1-1/U$, \item
$\Psi_{U}^{(j)} (t) \ll U^j$ for $j = 1,2, ...,K$.
\end{enumerate}
This weight function satisfies the requirements in Goldston and
Gonek [\ref{GG}]. One more thing to note is that
$$Re \hat{\Psi}_{U}(y) = \frac{\sin{2\pi y}}{2\pi y} \Bigl(\frac{\sin{2\pi \Delta
y}}{2\pi \Delta y}\Bigr)^{K+1}$$ where $\hat{f}(y) =
\int_{-\infty}^{\infty} f(t) e(yt) dt$.

We also need to study
$$S_{\alpha}^h(y) := \sum_{k \leq y} {\mathfrak S}(k) k^{\alpha}
\cos{(h\log{\frac{kx}{y}})} - \int_{0}^{y} u^\alpha
\cos{(h\log{\frac{ux}{y}})} du \mbox{ for } \alpha \geq 0,$$ and
$$T_{\alpha}^h(y) := \sum_{k >y} \frac{{\mathfrak S}(k)}{k^\alpha}
\cos{(h\log{\frac{kx}{y}})} - \int_{y}^{\infty} \frac{1}{u^\alpha}
\cos{(h\log{\frac{ux}{y}})} du \mbox{ for } \alpha>1.$$ Then from
[\ref{FG}],
\begin{equation}
\label{4.1} S_0(y) := S_0^0(y) = -\frac{1}{2}\log{y} +
O((\log{y})^{2/3}) = -\frac{1}{2}\log{y} + \epsilon(y).
\end{equation}
By partial summation and Lemma \ref{lemma1.7},
\begin{equation}
\label{4.2}
\begin{split}
S_\alpha^h(y) =& \epsilon(y) y^\alpha \cos{(h\log{x})} -
\frac{\alpha \cos{(h\log{x})}}{2(\alpha^2+h^2)} y^\alpha -
\frac{h \sin{(h\log{x})}}{2(\alpha^2+h^2)} y^\alpha \\
&- \int_{0}^{y} \epsilon(u) u^{\alpha-1} \Bigl[\alpha
\cos{(h\log{\frac{ux}{y}})} - h \sin{(h\log{\frac{ux}{y}})} \Bigr]
du,
\end{split}
\end{equation}
and
\begin{equation}
\label{4.3}
\begin{split}
T_\alpha^h(y) =& -\frac{\epsilon(y)}{y^\alpha} \cos{(h\log{x})} -
\frac{\alpha \cos{(h\log{x})}}{2(\alpha^2+h^2)} \frac{1}{y^\alpha}
+ \frac{h \sin{(h\log{x})}}{2(\alpha^2+h^2)} \frac{1}{y^\alpha} \\
&+ \int_{y}^{\infty} \frac{\epsilon(u)}{u^{\alpha+1}} \Bigl[\alpha
\cos{(h\log{\frac{ux}{y}})} + h \sin{(h\log{\frac{ux}{y}})} \Bigr]
du.
\end{split}
\end{equation}
Let
$$f(y) := \int_{0}^{y} \epsilon(u) - \frac{B}{2}\, du$$
where $B = -C_0 - \log{2\pi}$ and $C_0$ is Euler's constant. Note
that
\begin{equation}
\label{ineq} f(y) \ll y^{1/2 + \epsilon}
\end{equation}
(see Lemma 2.2 of [\ref{C3}]). From (\ref{4.2}) and (\ref{4.3}),
\begin{equation}
\label{4.8}
\begin{split}
& S_2^h(y) \frac{1}{y^3} + T_2^h(y) y \\
=& - \frac{2 \cos{(h\log{x})}}{4+h^2} \frac{1}{y} - \frac{1}{y^3}
\int_{0}^{y} u\epsilon(u) \Bigl[2\cos{(h\log{\frac{ux}{y}})} -
h\sin{(h\log{\frac{ux}{y}})} \Bigr] du \\
&+ y \int_{y}^{\infty} \frac{\epsilon(u)}{u^3}
\Bigl[2\cos{(h\log{\frac{ux}{y}})} + h\sin{(h\log{\frac{ux}{y}})}
\Bigr] du.
\end{split}
\end{equation}
\begin{lemma}
\label{lemma4.1}
\begin{equation*}
\begin{split}
I + J =&- \frac{1}{y^3} \int_{0}^{y} u\epsilon(u) \Bigl[
2\cos{(h\log{\frac{ux}{y}})} - h\sin{(h\log{\frac{ux}{y}})}
\Bigr] du \\
&+ y \int_{y}^{\infty} \frac{\epsilon(u)}{u^3}
\Bigl[2\cos{(h\log{\frac{ux}{y}})} + h\sin{(h\log{\frac{ux}{y}})}
\Bigr] du \\
=& - \frac{4f(y)}{y} \cos{(h\log{x})} \\
&+ \frac{1}{y^3} \int_{0}^{y} f(u) \Bigl[ (2-h^2)
\cos{(h\log{\frac{ux}{y}})} -
3h\sin{(h\log{\frac{ux}{y}})} \Bigr] du \\
&+ y \int_{y}^{\infty} \frac{f(u)}{u^4} \Bigl[ (6-h^2)
\cos{(h\log{\frac{ux}{y}})} + 5h\sin{(h\log{\frac{ux}{y}})} \Bigr]
du.
\end{split}
\end{equation*}
\end{lemma}

Proof: $I$ can be rewritten as
\begin{equation*}
\begin{split}
& -\frac{1}{y^3} \int_{0}^{y} u \Bigl(\epsilon(u) -
\frac{B}{2}\Bigr) \Bigl[2\cos{(h\log{\frac{ux}{y}})} -
h\sin{(h\log{\frac{ux}{y}})} \Bigr] du \\
-& \frac{B}{2} \frac{1}{y^3} \int_{0}^{y} u
\Bigl[2\cos{(h\log{\frac{ux}{y}})} - h\sin{(h\log{\frac{ux}{y}})}
\Bigr] du = - I_1 - I_2.
\end{split}
\end{equation*}
By a substitution $\Dst{v = \log{\frac{ux}{y}}}$ and Lemma
\ref{lemma1.7},
\begin{equation}
\label{4.4} I_2 = \frac{B}{2} \frac{1}{y} \cos{(h\log{x})}.
\end{equation}
By integration by parts and (\ref{ineq}),
\begin{equation}
\label{4.5}
\begin{split}
I_1 =& \frac{f(y)}{y^2} [2\cos{(h\log{x})} - h\sin{(h\log{x})}]
\\
&- \frac{1}{y^3} \int_{0}^{y} f(u) \Bigl[ (2-h^2)
\cos{(h\log{\frac{ux}{y}})} - 3h\sin{(h\log{\frac{ux}{y}})} \Bigr]
du.
\end{split}
\end{equation}
Similarly, $J$ can be rewritten as
\begin{equation*}
\begin{split}
& y \int_{y}^{\infty} \frac{\epsilon(u) - \frac{B}{2}}{u^3}
\Bigl[2\cos{(h\log{\frac{ux}{y}})} + h\sin{(h\log{\frac{ux}{y}})}
\Bigr] du \\
+& \frac{B}{2} y \int_{y}^{\infty} \frac{1}{u^3}
\Bigl[2\cos{(h\log{\frac{ux}{y}})} + h\sin{(h\log{\frac{ux}{y}})}
\Bigr] du = J_1 + J_2
\end{split}
\end{equation*}
By a substitution $\Dst{v = \log{\frac{ux}{y}}}$ and Lemma
\ref{lemma1.7},
\begin{equation}
\label{4.6} J_2 = \frac{B}{2} \frac{1}{y} \cos{(h\log{x})}.
\end{equation}
By integration by parts and (\ref{ineq}),
\begin{equation}
\label{4.7}
\begin{split}
J_1 =& - \frac{f(y)}{y^2} [2\cos{(h\log{x})} + h\sin{(h\log{x})}]
\\
&+ y \int_{y}^{\infty} \frac{f(u)}{u^4} \Bigl[ (6-h^2)
\cos{(h\log{\frac{ux}{y}})} + 5h\sin{(h\log{\frac{ux}{y}})} \Bigr]
du.
\end{split}
\end{equation}
(\ref{4.4}), (\ref{4.5}), (\ref{4.6}) and (\ref{4.7}) together
gives the lemma.
\begin{lemma}
\label{lemma4.2}
$$S_2^h(y) \frac{1}{y^3} + T_2^h(y) y = - \frac{2 \cos{(h\log{x})}}
{4+h^2} \frac{1}{y} - \frac{4f(y)}{y^2} \cos{(h\log{x})} + G_1(y)
+ G_2(y)$$ where
$$G_1(y) = \frac{1}{y^3} \int_{0}^{y} f(u) \Bigl[ (2-h^2)
\cos{(h\log{\frac{ux}{y}})} - 3h\sin{(h\log{\frac{ux}{y}})} \Bigr]
du, $$ and
$$G_2(y) = y \int_{y}^{\infty} \frac{f(u)}{u^4} \Bigl[ (6-h^2)
\cos{(h\log{\frac{ux}{y}})} + 5h\sin{(h\log{\frac{ux}{y}})} \Bigr]
du.$$
\end{lemma}

Proof: Combine (\ref{4.8}) and Lemma \ref{lemma4.1}.
\begin{lemma}
\label{lemma4.3} For any integer $n \geq 1$,
$$\int_{1}^{\infty} \frac{1}{y^n} Re \hat{\Psi}_U \Bigl(\frac{Ty}{2\pi x}\Bigr)
dy = \int_{1}^{\infty} \frac{1}{y^n}
\frac{\sin{\frac{T}{x}y}}{\frac{T}{x}y} dy + O\Bigl(\Delta
\log{\frac{1}{\Delta}}\Bigr).$$ When $n \neq 2$, the error term
can be replaced by $O(\Delta)$.
\end{lemma}

Proof: This is Lemma 3.3 in [\ref{C3}].
\begin{lemma}
\label{lemma4.4} If $F(y) \ll y^{-{3/2} + \epsilon}$ for $y \geq
1$, then
$$\int_{1}^{\infty} F(y) Re \hat{\Psi}_U \Bigl(\frac{Ty}{2 \pi x}\Bigr) dy =
\int_{1}^{\infty} F(y) \frac{\sin{\frac{T}{x}y}}{\frac{T}{x}y} dy
+ O(\Delta).$$
\end{lemma}

Proof: This is Lemma 3.4 in [\ref{C3}].
\section{Proof of Theorem \ref{theorem2}}
Throughout this section, we assume $\tau = T^{1-\epsilon} \leq T
/\log^M{T} \leq x \leq T^{2-2\epsilon}$, $U = \log^M{T}$ for
$M>2$, $H^{*}=\tau^{-2} x^{2/(1- \epsilon)}$, and $\Psi_{U}(t)$ is
defined as in the previous section. The implicit constants in the
error terms may depend on $\epsilon$, $K$ and $M$.

Our method is that of Goldston and Gonek [\ref{GG}] and it is very
similar to [\ref{C3}]. Let $s = \sigma + it$,
$$A_h(s) := \sum_{n \leq x} \frac{\Lambda(n) (1+n^{ih})}{n^s} \mbox{ and }
A_{h}^{*}(s) := \sum_{n>x} \frac{\Lambda(n) (1+n^{ih})}{n^s},$$
$$A(s) := \frac{1}{2} A_0(s) \mbox{ and } A^{*}(s) :=
\frac{1}{2} A_0^{*}(s).$$ By Lemma \ref{lemma1.4}, with slight
modifications, one has
\begin{equation*}
\begin{split}
4\pi F_h(x,T) =& \int_{0}^{T} \Big| \frac{1}{x} \Bigl(
A_h(-\frac{1}{2} +it) - \int_{1}^{x} (1+u^{ih}) u^{1/2-it} du
\Bigr) \\
&+ x \Bigl( A_h^{*}(\frac{3}{2}+it) - \int_{x}^{\infty} (1+u^{ih})
u^{-3/2-it} du \Bigr) \Big|^2 dt \\
-& 2 \int_{0}^{T} \Big| \frac{1}{x} \Bigl( A(-\frac{1}{2} +it) -
\int_{1}^{x} u^{1/2-it} du
\Bigr) \\
&+ x \Bigl( A^{*}(\frac{3}{2}+it) - \int_{x}^{\infty} u^{-3/2-it}
du \Bigr) \Big|^2 dt + O\bigl(\tilde{h}\log^3{T}\bigr).
\end{split}
\end{equation*}
Inserting $\Psi_{U} (t/T)$ into the integral and extending the
range of integration to the whole real line, we have
\begin{equation}
\label{5.1}
\begin{split}
4\pi F(x,T) =& \frac{1}{x^2} I_1 (x, T) + x^2 I_2 (x, T) -
\frac{2}{x^2} I_3(x, T) - 2x^2 I_4(x, T) \\
&+ O \Bigl(\frac{T(\log{T})^2}{U}\Bigr) + O \Bigl(
\frac{x^{1+6\epsilon}}{T}\Bigr)
\end{split}
\end{equation}
where
\begin{equation*}
\begin{split}
I_1(x,T) =& \int_{-\infty}^{\infty} \Big| A_h(-\frac{1}{2}+it) -
\int_{1}^{x} (1+u^{ih}) u^{1/2-it} du \Big|^2
\Psi_{U}\Bigl(\frac{t}{T} \Bigr) dt, \\
I_2(x,T) =& \int_{-\infty}^{\infty} \Big| A_h^{*}(\frac{3}{2}+it)
- \int_{x}^{\infty} (1+u^{ih}) u^{-3/2-it} du \Big|^2
\Psi_{U}\Bigl(\frac{t}{T} \Bigr) dt, \\
I_3(x,T) =& \int_{-\infty}^{\infty} \Big| A(-\frac{1}{2}+it) -
\int_{1}^{x} u^{1/2-it} du \Big|^2 \Psi_{U}\Bigl(\frac{t}{T}
\Bigr) dt, \\
I_4(x,T) =& \int_{-\infty}^{\infty} \Big| A^{*}(\frac{3}{2}+it) -
\int_{x}^{\infty} u^{-3/2-it} du \Big|^2 \Psi_{U}\Bigl(\frac{t}{T}
\Bigr) dt
\end{split}
\end{equation*}
by Lemma $1$ of [\ref{GGOS}] with modification $V = -T/U$ and $T-T
/U$, and $W = 2T/U$. The contribution from the cross terms are
estimated via Theorem $3$ of [\ref{GG}]. Note that by partial
summation with the Riemann Hypothesis and {\bf TPC},
\begin{equation*}
\begin{split}
& \sum_{n \leq x} \Lambda(n) (1+n^{ih}) = \int_{1}^{x}
(1+u^{ih})\,
du + O(\tilde{h} x^{1/2 + \epsilon}), \\
& \sum_{n \leq x} \Lambda(n) \Lambda(n+k) (1+n^{ih})
(1+(n+k)^{-ih}) \\
=& {\mathfrak S}(k) \int_{1}^{x} (1+u^{ih})(1+(u+k))^{-ih}\, du +
O(\tilde{h} x^{1/2 + \epsilon}).
\end{split}
\end{equation*}
By Corollary $1$ of [\ref{GG}] (see also the calculations at the
end of [\ref{GG}] and [\ref{GGOS}]),
\begin{equation*}
\begin{split}
I_1(x,T) =& \hat{\Psi}_{U}(0) T \sum_{n \leq x} \Lambda^2 (n) n |1+n^{ih}|^2 \\
&+ 4 \pi \Bigl(\frac{T}{2 \pi}\Bigr)^3 \int_{T/2\pi x}^{\infty}
\Bigl[ \sum_{k \leq 2 \pi x v /T} {\mathfrak S} (k) k^2 \bigl(1 +
(\frac{kT}{2\pi v})^{ih}\bigr) \\
&\times \bigl(1 + (\frac{kT}{2\pi v} + k)^{-ih}\bigr) \Bigr] Re
\hat{\Psi}_{U}(v) \frac{dv}{v^3} \\
&- 4 \pi \Bigl(\frac{T}{2 \pi}\Bigr)^3 \int_{T/2\pi \tau
x}^{\infty} \Bigl[ \int_{0}^{2\pi x v/T} u^2 \big|1 + (\frac{uT}{2
\pi v})^{ih}
\big|^2 du \Bigr] Re \hat{\Psi}_{U}(v) \frac{dv}{v^3} \\
&+ O\Bigl(\tilde{h}\frac{x^{3+6\epsilon}}{T}\Bigr) + O(\tilde{h}
x^{{5/2}+7\epsilon}).
\end{split}
\end{equation*}
Note that
\begin{equation}
\label{5.2}
\begin{split}
& \Bigl(1 + (\frac{kT}{2\pi v})^{ih}\Bigr) \Bigl(1 +
(\frac{kT}{2\pi v} + k)^{-ih}\Bigr) \\
=& \Bigl|1 + (\frac{kT}{2\pi v})^{ih}\Bigr|^2 + \Bigl(1 +
(\frac{kT}{2\pi v})^{ih}\Bigr) \Bigl((\frac{kT}{2\pi v} + k)^{-ih}
- (\frac{kT}{2\pi v})^{-ih}\Bigr) \\
=& \Bigl|1 + (\frac{kT}{2\pi v})^{ih}\Bigr|^2 + \Bigl(1 +
(\frac{kT}{2\pi v})^{ih}\Bigr) \Bigl(\frac{kT}{2\pi v}\Bigr)^{-ih}
\Bigl((1 + \frac{2\pi v}{T})^{-ih} - 1 \Bigr) \\
=& \Bigl|1 + (\frac{kT}{2\pi v})^{ih}\Bigr|^2 +
O\Bigl(\min{\bigl(\frac{h v}{T}, 1\bigr)}\Bigr).
\end{split}
\end{equation}
Thus,
\begin{equation*}
\begin{split}
& \int_{T/2\pi x}^{\infty} \Bigl[ \sum_{k \leq 2 \pi x v /T}
{\mathfrak S} (k) k^2 \bigl(1 + (\frac{kT}{2\pi v})^{ih}\bigr)
\bigl(1 + (\frac{kT}{2\pi v} + k)^{-ih}\bigr) \Bigr] Re
\hat{\Psi}_{U}(v) \frac{dv}{v^3} \\
=& \int_{T/2\pi x}^{\infty} \Bigl[ \sum_{k \leq 2 \pi x v /T}
{\mathfrak S} (k) k^2 \bigl|1 + (\frac{kT}{2\pi v})^{ih}\bigr|^2
\Bigr] Re \hat{\Psi}_{U}(v) \frac{dv}{v^3} \\
&+ O\Bigl(\int_{T/2\pi x}^{T^{2-\epsilon}/x} \bigl( \frac{xv}{T}
\bigr)^3 \frac{hv}{T} \frac{1}{v} \frac{dv}{v^3} +
\int_{T^{2-\epsilon}/x}^{\infty} \bigl( \frac{xv}{T} \bigr)^3
\frac{1}{\Delta v^2} \frac{dv}{v^3} \Bigr) \\
=& \int_{T/2\pi x}^{\infty} \Bigl[ \sum_{k \leq 2 \pi x v /T}
{\mathfrak S} (k) k^2 \bigl|1 + (\frac{kT}{2\pi v})^{ih}\bigr|^2
\Bigr] Re \hat{\Psi}_{U}(v) \frac{dv}{v^3} + O\Bigl(\frac{h
x^2}{T^{2+\epsilon}} + \frac{x^4}{\Delta T^{5-\epsilon}}\Bigr)
\end{split}
\end{equation*}
as $\Dst{\sum_{k \leq x} {\mathfrak S}(k) \sim x}$ and $\Dst{Re
\hat{\Psi}_{U}(v) \ll \min{(\frac{1}{v}, \frac{1}{\Delta v^2})}}$.
Therefore,
\begin{equation*}
\begin{split}
I_1(x,T) =& T \sum_{n \leq x} \Lambda^2 (n) n (2 + 2\cos{(h\log{n})}) \\
&+ 4 \pi \Bigl(\frac{T}{2 \pi}\Bigr)^3 \int_{T/2\pi x}^{\infty}
\Bigl[ \sum_{k \leq 2 \pi x v /T} {\mathfrak S} (k) k^2 \bigl(2 +
2\cos{(h \log{\frac{kT}{2\pi v}})}\bigr) \\
&- \int_{0}^{2\pi x v/T} u^2 \bigl(2 + 2\cos{(h
\log{\frac{uT}{2\pi v}})}\big) du \Bigr] Re \hat{\Psi}_{U}(v)
\frac{dv}{v^3} \\
&- 4 \pi \Bigl(\frac{T}{2 \pi}\Bigr)^3 \int_{0}^{T/2\pi x}
\int_{0}^{2\pi x v/T} u^2 \bigl(2 + 2\cos{(h \log{\frac{uT}{2\pi
v}})} \big) du Re \hat{\Psi}_{U}(v) \frac{dv}{v^3} \\
&+ O\Bigl(\frac{\tilde{h} x^{3+6\epsilon}}{T}\Bigr) + O(\tilde{h}
x^{{5/2}+7\epsilon}) + O(\tilde{h} x^2 T^{1-\epsilon}) +
O\Bigl(\frac{x^4}{\Delta T^{2-\epsilon}}\Bigr).
\end{split}
\end{equation*}
Similarly, by Corollary $2$ of [\ref{GG}],
\begin{equation*}
\begin{split}
I_2(x,T) =& T \sum_{x<n} \frac{\Lambda^2(n)}{n^3} (2 + 2\cos{(h\log{n})}) \\
&+ \frac{8 \pi^2}{T} \int_{0}^{T H^{*}/2 \pi x} \Bigl[\sum_{2\pi x
v/T \leq k \leq H^{*}} \frac{{\mathfrak S}(k)}{k^2} \bigl(2 +
2\cos{(h \log{\frac{kT}{2\pi v}})}\bigr) \\
&- \int_{2\pi xv/T}^{H^{*}} \frac{1}{u^2} \bigl(2 + 2\cos{(h
\log{\frac{uT}{2\pi v}})}\bigr) du \Bigr] Re \hat{\Psi}_{U} (v) v dv \\
&+ O(\tilde{h} T^{-1} x^{-1+6\epsilon}) + O(\tilde{h} x^{-{3/2} +
7\epsilon}) + O(\tilde{h} T^{1-{\epsilon/2} } x^{-2}) +
O\Bigl(\frac{\tilde{h} H^*}{\Delta x^2} \Bigr)
\end{split}
\end{equation*}
where the last error term comes from the error term in
(\ref{5.2}). $I_3(x,T)$ and $I_4(x,T)$ are computed in [\ref{C3}]
or one can simply set $h=0$ in $I_1(x,T)$ and $I_2(x,T)$, and
divide by $4$. Putting these into (\ref{5.1}) with a substitution
$\Dst{y = \frac{2\pi x v}{T}}$ and using Lemma \ref{lemma1.8},
\begin{equation*}
\begin{split}
& 4\pi F_h(x,T) \\
&= 2T \Bigl[\frac{4 \cos{(h\log{x})}}{4 + h^2} \log{x}
- \frac{8h \sin{(h \log{x})}}{(4 + h^2)^2} \Bigr] \\
&+ 4T \int_{1}^{\infty} \Bigl[ \sum_{k \leq y} {\mathfrak S} (k)
k^2 \cos{(h \log{\frac{kx}{y}})} - \int_{0}^{y} u^2 \cos{(h
\log{\frac{ux}{y}})} du \Bigr]
Re \hat{\Psi}_{U}\Bigl(\frac{Ty}{2\pi x}\Bigr) \frac{dy}{y^3} \\
&- 4T \int_{0}^{1} \int_{0}^{y} u^2 \cos{(h \log{\frac{ux}{y}})} du
Re \hat{\Psi}_{U}\Bigl(\frac{Ty}{2\pi x}\Bigr) \frac{dy}{y^3} \\
&+ 4T \int_{1}^{H^{*}} \Bigl[\sum_{y \leq k \leq H^{*}}
\frac{{\mathfrak S}(k) \cos{(h \log{\frac{kx}{y}})}}{k^2} -
\int_{y}^{H^{*}} \frac{\cos{(h \log{\frac{ux}{y}})}}{u^2} du
\Bigr] Re \hat{\Psi}_{U}\Bigl(\frac{Ty}{2\pi x}\Bigr) y dy \\
&+ 4T \int_{0}^{1} \Bigl[\sum_{k \leq H^{*}} \frac{{\mathfrak
S}(k) \cos{(h \log{\frac{kx}{y}})}}{k^2} - \int_{y}^{H^{*}}
\frac{\cos{(h \log{\frac{ux}{y}})}}{u^2} du \Bigr]
Re \hat{\Psi}_{U}\Bigl(\frac{Ty}{2\pi x}\Bigr) y dy \\
&+ O\Bigl(\frac{\tilde{h} x^{1+6\epsilon}}{T}\Bigr) + O(\tilde{h}
x^{1/2 + 7\epsilon}) + O\Bigl(\frac{\tilde{h}
x^2}{T^{2-2\epsilon}}\Bigr) +
O\Bigl(\frac{\tilde{h}T}{\log^{M-2}{T}}\Bigr)
\end{split}
\end{equation*}
From Lemma \ref{lemma1.7},
\begin{equation}
\label{5.4} \int e^{ax}\cos{bx} = \frac{a}{a^2+b^2} e^{ax}\cos{bx}
+ \frac{b}{a^2+b^2} e^{ax}\sin{bx}.
\end{equation}
Also,
\begin{equation}
\label{5.5}
\begin{split}
\int_{0}^{1} Re \hat{\Psi}_{U}\Bigl(\frac{Ty}{2\pi x}\Bigr) dy =&
\frac{x}{T} \int_{0}^{T/x} \frac{\sin{u}}{u} (1 + O(\Delta^2 u^2)) du \\
=& \frac{x}{T} \int_{0}^{T/x} \frac{\sin{u}}{u} du +
O\Bigl(\frac{\Delta^2 T}{x}\Bigr).
\end{split}
\end{equation}
By appropriate change of variables, (\ref{5.4}) and (\ref{5.5}),
\begin{equation*}
\begin{split}
& \int_{0}^{1} \int_{0}^{y} u^2 \cos{(h \log{\frac{ux}{y}})} du
Re \hat{\Psi}_{U}\Bigl(\frac{Ty}{2\pi x}\Bigr) \frac{dy}{y^3} \\
=& \frac{x}{T} \int_{0}^{T/x} \frac{\sin{u}}{u} du
\Bigl[\frac{3}{9+h^2}\cos{(h\log{x})} +
\frac{h}{9+h^2}\sin{(h\log{x})}\Bigr] + O\Bigl(\frac{\Delta^2
T}{x}\Bigr). \\
& \int_{0}^{1} \Bigl[\sum_{k \leq H^{*}} \frac{{\mathfrak S}(k)
\cos{(h \log{\frac{kx}{y}})}}{k^2} - \int_{y}^{H^{*}}
\frac{\cos{(h \log{\frac{ux}{y}})}}{u^2} du \Bigr] Re
\hat{\Psi}_{U}\Bigl(\frac{Ty}{2\pi x}\Bigr) y dy \\
=& \sum_{k=1}^{\infty} \frac{{\mathfrak S}(k)}{k^2} \int_{0}^{1} y
\cos{(h\log{\frac{kx}{y}})}
\frac{\sin{\frac{Ty}{x}}}{\frac{Ty}{x}} dy +
O\Bigl(\frac{1}{H^*}\Bigr) + O\Bigl(\frac{\Delta^2 T}{x}\Bigr) \\
&- \frac{x}{T} \int_{0}^{T/x} \frac{\sin{u}}{u} du
\Bigl[\frac{1}{1+h^2}\cos{(h\log{x})} -
\frac{h}{1+h^2}\sin{(h\log{x})}\Bigr] + O\Bigl(\frac{\Delta^2
T}{x}\Bigr).
\end{split}
\end{equation*}
Therefore, with the notation $\Dst{S_\alpha^h(y)}$ and
$\Dst{T_\alpha^h(y)}$,
\begin{equation*}
\begin{split}
& 4\pi F_h(x,T) \\
&= 2T \Bigl[\frac{4 \cos{(h\log{x})}}{4 + h^2} \log{x}
- \frac{8h \sin{(h \log{x})}}{(4 + h^2)^2} \Bigr] \\
&+ 4T \int_{1}^{\infty} S_2^h(y)
Re \hat{\Psi}_{U}\Bigl(\frac{Ty}{2\pi x}\Bigr) \frac{dy}{y^3} \\
&+ 4T \int_{1}^{H^{*}} \bigl(T_2^h(y) - T_2^h(H^*)\bigr)
Re \hat{\Psi}_{U}\Bigl(\frac{Ty}{2\pi x}\Bigr) y dy \\
&- 4x \int_{0}^{T/x} \frac{\sin{u}}{u} du
\Bigl[\frac{3\cos{(h\log{x})}}{9+h^2} +
\frac{h\sin{(h\log{x})}}{9+h^2}\Bigr] \\
&- 4x \int_{0}^{T/x} \frac{\sin{u}}{u} du
\Bigl[\frac{\cos{(h\log{x})}}{1+h^2} -
\frac{h\sin{(h\log{x})}}{1+h^2}\Bigr] \\
&+ 4T \sum_{k=1}^{\infty} \frac{{\mathfrak S}(k)}{k^2}
\int_{0}^{1} y \cos{(h\log{\frac{kx}{y}})}
\frac{\sin{\frac{Ty}{x}}} {\frac{Ty}{x}} dy \\
&+ O\Bigl(\frac{\tilde{h} x^{1+6\epsilon}}{T}\Bigr) + O(\tilde{h}
x^{1/2 + 7\epsilon}) + O\Bigl(\frac{\tilde{h}
x^2}{T^{2-2\epsilon}}\Bigr) +
O\Bigl(\frac{\tilde{h}T}{\log^{M-2}{T}}\Bigr).
\end{split}
\end{equation*}
By (\ref{4.1}) and (\ref{4.3}), $\Dst{T_2^h(H^*) \ll \frac{h
(\log{H^*})^{2/3}}{(H^*)^2}}$. It follows that the contribution
from $T_2^h(H^*)$ in the second integral is $O(h T^{-\epsilon})$.
Also, one can extend the upper limit of the second integral to
$\infty$ with an error $O(h T^{-\epsilon})$ by (\ref{4.1}) and
(\ref{4.3}) again. Finally, we obtain the theorem by applying
Lemma \ref{lemma4.2}, Lemma \ref{lemma4.3}, Lemma \ref{lemma4.4},
(\ref{ineq}) and dividing by $4\pi$.
\section{Proof of Theorem \ref{theorem3} and \ref{theorem4}}

Proof of Theorem \ref{theorem3}: It follows directly from Theorem
\ref{theorem2} by observing that all the other main terms besides
the first one are $O(x)$ because of (\ref{ineq}).

\bigskip

Before proving Theorem \ref{theorem4}, we need the following
lemmas.
\begin{lemma}
\label{lemma6.1}
$$\int_{1}^{\infty} \frac{\sin{ax}}{x^{2n}} dx = \frac{a^{2n-1}}{(2n-1)!}
\Bigl[\sum_{k=1}^{2n-1} \frac{(2n-k-1)!}{a^{2n-k}}
\sin{\bigl(a+(k-1)\frac{\pi}{2}\bigr)} + (-1)^n ci(a) \Bigr]$$
where $ci(x) = -\int_{x}^{\infty} \frac{\cos{t}}{t}dt = C_0 +
\log{x} + \int_{0}^{x} \frac{\cos{t}-1}{t} dt$ and $C_0$ is
Euler's constant.
\end{lemma}

Proof: This is formula $3.761(3)$ on P.430 of [\ref{GR}] which can
be proved by integration by parts repeatedly.

\begin{lemma}
\label{lemma6.2} If $F(y) \ll y^{-{3/2}+\epsilon}$ for $y \geq 1$,
then for $T \leq x$,
$$\int_{1}^{\infty} F(y) \frac{\sin{\frac{T}{x}y}}{\frac{T}{x}y} dy =
\int_{1}^{\infty} F(y) dy +
O\Bigl(\bigl(\frac{T}{x}\bigr)^{{1/2}-\epsilon} \Bigr).$$
\end{lemma}

Proof: This is Lemma 5.2 in [\ref{C3}].

\begin{lemma}
\label{lemma6.3}
\begin{equation*}
\begin{split}
I =& \int_{1}^{\infty} \frac{1}{y^3} \int_{0}^{y} f(u)
\Bigl[(2-h^2) \cos{(h\log{\frac{ux}{y}})} - 3h
\sin{(h\log{\frac{ux}{y}})} \Bigr] du\, dy \\
=& \int_{0}^{1} f(u) [\cos{(h\log{ux})} - h\sin{(h\log{ux})}]\, du
\\
&+ \int_{1}^{\infty} \frac{f(u)}{u^2}\, du [\cos{(h\log{x})} -
h\sin{(h\log{x})}].
\end{split}
\end{equation*}
\end{lemma}

Proof: Because of (\ref{ineq}), we can change the order of
integration.
\begin{equation*}
\begin{split}
I =& \int_{0}^{1} f(u) \int_{1}^{\infty} \frac{1}{y^3}
\Bigl[(2-h^2) \cos{(h\log{\frac{ux}{y}})} - 3h
\sin{(h\log{\frac{ux}{y}})} \Bigr] dy\, du \\
&+ \int_{1}^{\infty} f(u) \int_{u}^{\infty} \frac{1}{y^3}
\Bigl[(2-h^2) \cos{(h\log{\frac{ux}{y}})} - 3h
\sin{(h\log{\frac{ux}{y}})} \Bigr] dy\, du \\
=& \int_{0}^{1} f(u) \Bigl\{ (2-h^2) \Bigl[ \frac{2}{4+h^2}
\cos{(h\log{ux})} + \frac{h}{4+h^2} \sin{(h\log{ux})} \Bigr] \\
&- 3h \Bigl[ \frac{2}{4+h^2} \sin{(h\log{ux})} - \frac{h}{4+h^2}
\cos{(h\log{ux})} \Bigr] \Bigr\} du \\
& + \int_{1}^{\infty} \frac{f(u)}{u^2} \Bigl\{ (2-h^2) \Bigl[
\frac{2}{4+h^2} \cos{(h\log{x})} + \frac{h}{4+h^2}
\sin{(h\log{x})} \Bigr] \\
&- 3h \Bigl[ \frac{2}{4+h^2} \sin{(h\log{x})} - \frac{h}{4+h^2}
\cos{(h\log{x})} \Bigr] \Bigr\} du.
\end{split}
\end{equation*}
by substituting $v = \log{\frac{ux}{y}}$ and applying Lemma
\ref{lemma1.7}. The lemma follows after some simple algebra.

\begin{lemma}
\label{lemma6.4}
\begin{equation*}
\begin{split}
J =& \int_{1}^{\infty} y \int_{y}^{\infty} \frac{f(u)}{u^4}
\Bigl[(6-h^2) \cos{(h\log{\frac{ux}{y}})} + 5h
\sin{(h\log{\frac{ux}{y}})} \Bigr] du\, dy \\
=& - \int_{1}^{\infty} \frac{f(u)}{u^4} [3\cos{(h\log{ux})} +
h\sin{(h\log{ux})}]\, du \\
&+ \int_{1}^{\infty} \frac{f(u)}{u^2}\, du [3\cos{(h\log{x})} +
h\sin{(h\log{x})}].
\end{split}
\end{equation*}
\end{lemma}

Proof: Again, because of (\ref{ineq}), we can change the order of
integration.
\begin{equation*}
\begin{split}
J =& \int_{1}^{\infty} \int_{1}^{u} y \Bigl[(6-h^2)
\cos{(h\log{\frac{ux}{y}})} + 5h \sin{(h\log{\frac{ux}{y}})}
\Bigr] dy\, du \\
=& \int_{1}^{\infty} \frac{f(u)}{u^4} \Bigl\{(6-h^2) \Bigl[
\frac{-2}{4+h^2} \cos{(h\log{ux})} + \frac{h}{4+h^2}
\sin{(h\log{ux})} \Bigr] \\
&+ 5h \Bigl[\frac{-2}{4+h^2} \sin{(h\log{ux})} - \frac{h}{4+h^2}
\cos{(h\log{ux})} \Bigr] \\
&- (6-h^2) \Bigl[ \frac{-2}{4+h^2} \cos{(h\log{x})} +
\frac{h}{4+h^2} \sin{(h\log{x})} \Bigr] \\
&- 5h \Bigl[\frac{-2}{4+h^2} \sin{(h\log{x})} - \frac{h}{4+h^2}
\cos{(h\log{x})} \Bigr] \Bigr\} du.
\end{split}
\end{equation*}
by substituting $v = \log{\frac{ux}{y}}$ and applying Lemma
\ref{lemma1.7}. The lemma follows after some simple algebra.

\begin{lemma}
\label{lemma6.5}
\begin{equation*}
\begin{split}
S =& \sum_{k=1}^{\infty} \frac{{\mathfrak S}(k)}{k^2} \int_{0}^{1}
y \cos{(h\log{\frac{kx}{y}})}\, dy \\
=& \Bigl[\frac{1}{1+h^2} \cos{(h\log{x})} - \frac{h}{1+h^2}
\sin{(h\log{x})} \Bigr] \\
&- \Bigl[\frac{4-h^2}{2(4+h^2)^2}
\cos{(h\log{x})} - \frac{2h}{(4+h^2)^2} \sin{(h\log{x})}\Bigr] \\
&+ \frac{B}{2} \Bigl[\frac{2}{4+h^2} \cos{(h\log{x})} -
\frac{h}{4+h^2} \sin{(h\log{x})}\Bigr] + \Bigl(1 +
\frac{B}{2}\Bigr) \cos{(h\log{x})} \\
&+ \int_{1}^{\infty} \frac{f(u)}{u^4} \Bigl[ 3 \cos{(h\log{ux})} +
h \sin{(h\log{ux})} \Bigr] du.
\end{split}
\end{equation*}
\end{lemma}

Proof: By substituting $v = \log{\frac{kx}{y}}$ and Lemma
\ref{lemma1.7},
$$S = \frac{2}{4+h^2} \sum_{k=1}^{\infty}
\frac{{\mathfrak S}(k)}{k^2} \cos{(h\log{kx})} - \frac{h}{4+h^2}
\sum_{k=1}^{\infty} \frac{{\mathfrak S}(k)}{k^2}
\sin{(h\log{kx})}.$$ Recall the definition of $S_0(u)$ from
(\ref{4.1}) and use partial summation,
\begin{equation*}
\begin{split}
S =& \frac{2}{4+h^2} \int_{1}^{\infty} \frac{S_0(u) + u}{u^3} [2
\cos{(h\log{ux})} + h \sin{(h\log{ux})}]\, du \\
&- \frac{h}{4+h^2} \int_{1}^{\infty} \frac{S_0(u) + u}{u^3} [-h
\cos{(h\log{ux})} + 2 \sin{(h\log{ux})}]\, du \\
=& \int_{1}^{\infty} \frac{S_0(u) + u}{u^3} \cos{(h\log{ux})}\, du
\\
=& \int_{1}^{\infty} \frac{u - \frac{1}{2}\log{u} + \epsilon(u)}
{u^3} \cos{(h\log{ux})}\, du \\
=& \int_{1}^{\infty} \frac{1}{u^2} \cos{(h\log{ux})}\, du -
\frac{1}{2} \int_{1}^{\infty} \frac{\log{u}}{u^3}
\cos{(h\log{ux})}\, du \\
& + \frac{B}{2} \int_{1}^{\infty} \frac{1}{u^3}
\cos{(h\log{ux})}\, du + \int_{1}^{\infty} \frac{\epsilon(u) -
\frac{B}{2}}{u^3} \cos{(h\log{ux})}\, du \\
=& I_1 - \frac{1}{2} I_2 + \frac{B}{2} I_3 + I_4.
\end{split}
\end{equation*}
By appropriate substitution and Lemma \ref{lemma1.7},
\begin{equation*}
\begin{split}
I_1 =& \frac{1}{1+h^2} \cos{(h\log{x})} - \frac{h}{1+h^2}
\sin{(h\log{x})}, \\
I_2 =& \frac{4-h^2}{(4+h^2)^2} \cos{(h\log{x})} -
\frac{2h}{(4+h^2)^2} \sin{(h\log{x})}, \\
I_3 =& \frac{2}{4+h^2} \cos{(h\log{x})} - \frac{h}{4+h^2}
\sin{(h\log{x})}.
\end{split}
\end{equation*}
Finally, by integration by parts,
\begin{equation*}
\begin{split}
I_4 =& \int_{1}^{\infty} \frac{\cos{(h\log{ux})}}{u^3} df(u) \\
=& \Bigl(1 + \frac{B}{2}\Bigr) \cos{(h\log{x})} +
\int_{1}^{\infty} \frac{f(u)}{u^4} [3\cos{(h\log{ux})} +
h\sin{(h\log{ux})}]\, du
\end{split}
\end{equation*}
because $f(1) = -1 - \frac{B}{2}$. Combining the results for
$I_1$, $I_2$, $I_3$ and $I_4$, we have the lemma.

\begin{lemma}
\label{lemma6.6}
\begin{equation*}
\begin{split}
& \int_{0}^{1} f(u) [\cos{(h\log{ux})} - h\sin{(h\log{ux})}]\, du \\
=& -\frac{1}{2} \Bigl[\frac{4+3h^2}{(4+h^2)^2} \cos{(h\log{x})} +
\frac{h^3}{(4+h^2)^2} \sin{(h\log{x})}\Bigr] \\
&- \Bigl(\frac{1}{2} + \frac{B}{2}\Bigr) \Bigl[\frac{2+h^2}{4+h^2}
\cos{(h\log{x})} - \frac{h}{4+h^2} \sin{(h\log{x})} \Bigr] \\
&- \frac{1}{2} \Bigl[\frac{3+h^2}{9+h^2} \cos{(h\log{x})} -
\frac{2h}{9+h^2} \sin{(h\log{x})}\Bigr].
\end{split}
\end{equation*}
\end{lemma}

Proof: The key is $\epsilon(u) = \frac{1}{2}\log{u} - u$ when $0
\leq u \leq 1$ (see (\ref{4.1})). So,
$$f(u) = \int_{0}^{u} \epsilon(v) - \frac{B}{2}\, dv =
\frac{1}{2}u\log{u} - \Bigl(\frac{1}{2} + \frac{B}{2}\Bigr)u -
\frac{1}{2}u^2.$$ Putting this into the integral and evaluating
the integral piece by piece with suitable substitution and Lemma
\ref{lemma1.7}, one gets the lemma.

\bigskip

Proof of Theorem \ref{theorem4}: First observe that when $T \leq x
\leq T^{2 - 29\epsilon}$, the error terms in Theorem
\ref{theorem2} is $O(\frac{h T}{\log^{M-2}{T}})$. Rewrite Theorem
\ref{theorem2} as
$$F_h(x,T) = T_1 + T_2 + T_3 + T_4 + T_5 + O\Bigl(\frac{\tilde{h} T}
{\log^{M-2}{T}}\Bigr).$$ Since $\frac{\sin{u}}{u} = 1 + O(u^2)$,
\begin{equation*}
\begin{split}
T_3 &= - \frac{T}{\pi} \Bigl[\frac{3\cos{(h\log{x})}}{9+h^2} +
\frac{h\sin{(h\log{x})}}{9+h^2}\Bigr] + O\Bigl(T
\bigl(\frac{T}{x}\bigr)^2\Bigr), \\
T_4 &= - \frac{T}{\pi} \Bigl[\frac{\cos{(h\log{x})}}{1+h^2} -
\frac{h\sin{(h\log{x})}}{1+h^2}\Bigr]+ O\Bigl(T
\bigl(\frac{T}{x}\bigr)^2\Bigr).
\end{split}
\end{equation*}
By Lemma \ref{lemma6.5},
\begin{equation*}
\begin{split}
T_5 =& \frac{T}{\pi} \sum_{k=1}^{\infty} \frac{{\mathfrak
S}(k)}{k^2} \int_{0}^{1} y \cos{(h\log{\frac{kx}{y}})}\, dy +
O\Bigl(T \bigl(\frac{T}{x}\bigr)^2\Bigr) \\
=& \frac{T}{\pi} \Bigl[\frac{\cos{(h\log{x})}}{1+h^2} -
\frac{h\sin{(h\log{x})}}{1+h^2} \Bigr] \\
&- \frac{T}{\pi} \Bigl[\frac{4-h^2}{2(4+h^2)^2}
\cos{(h\log{x})} - \frac{2h}{(4+h^2)^2} \sin{(h\log{x})}\Bigr] \\
&+ \frac{T}{\pi} \frac{B}{2} \Bigl[\frac{2\cos{(h\log{x})}}{4+h^2}
- \frac{h\sin{(h\log{x})}}{4+h^2} \Bigr] +
\frac{T}{\pi} \Bigl(1 + \frac{B}{2}\Bigr) \cos{(h\log{x})} \\
&+ \frac{T}{\pi} \int_{1}^{\infty} \frac{f(u)}{u^4} \Bigl[ 3
\cos{(h\log{ux})} + h \sin{(h\log{ux})} \Bigr] du + O\Bigl(T
\bigl(\frac{T}{x}\bigr)^2\Bigr).
\end{split}
\end{equation*}
By Lemma \ref{lemma4.3}, (\ref{ineq}), Lemma \ref{lemma6.2}, Lemma
\ref{lemma6.3}, Lemma \ref{lemma6.4} and Lemma \ref{lemma6.6},
\begin{equation*}
\begin{split}
T_2 =& -\frac{2x\cos{(h\log{x})}}{\pi(4+h^2)} \int_{1}^{\infty}
\frac{\sin{\frac{Ty}{x}}}{y^2} dy - \frac{4T}{\pi}
\cos{(h\log{x})} \int_{1}^{\infty} \frac{f(y)}{y^2} dy \\
&+ \frac{T}{\pi} \int_{1}^{\infty} \bigl(G_1(y) + G_2(y)\bigr)\,
dy + O\Bigl(T \bigl(\frac{T}{x}\bigr)^{1/2-\epsilon}\Bigr) +
O\Bigl(\frac{\tilde{h}T}{\log^M{T}}\Bigr) \\
=& -\frac{2x\cos{(h\log{x})}}{\pi(4+h^2)} \int_{1}^{\infty}
\frac{\sin{\frac{Ty}{x}}}{y^2} dy + O\Bigl(\frac{\tilde{h}T}
{\log^M{T}}\Bigr)\\
&+ \frac{T}{\pi} \int_{0}^{1} f(u) [\cos{(h\log{ux})} -
h\sin{(h\log{ux})}]\, du \\
&- \frac{T}{\pi} \int_{1}^{\infty} \frac{f(u)}{u^4}
[3\cos{(h\log{ux})} + h\sin{(h\log{ux})}]\, du +
O\Bigl(T \bigl(\frac{T}{x}\bigr)^{1/2-\epsilon}\Bigr)\\
=& -\frac{2x\cos{(h\log{x})}}{\pi(4+h^2)} \int_{1}^{\infty}
\frac{\sin{\frac{Ty}{x}}}{y^2} dy + O\Bigl(\frac{\tilde{h}T}
{\log^M{T}}\Bigr)\\
&- \frac{T}{\pi} \frac{1}{2} \Bigl[\frac{4+3h^2}{(4+h^2)^2}
\cos{(h\log{x})} + \frac{h^3}{(4+h^2)^2} \sin{(h\log{x})}
\Bigr] \\
&- \frac{T}{\pi} \Bigl(\frac{1}{2} + \frac{B}{2}\Bigr)
\Bigl[\frac{2+h^2}{4+h^2}
\cos{(h\log{x})} - \frac{h}{4+h^2} \sin{(h\log{x})} \Bigr] \\
&- \frac{T}{\pi} \frac{1}{2} \Bigl[\frac{3+h^2}{9+h^2}
\cos{(h\log{x})} - \frac{2h}{9+h^2} \sin{(h\log{x})}\Bigr] \\
&- \frac{T}{\pi} \int_{1}^{\infty} \frac{f(u)}{u^4}
[3\cos{(h\log{ux})} + h\sin{(h\log{ux})}]\, du + O\Bigl(T
\bigl(\frac{T}{x}\bigr)^{1/2-\epsilon}\Bigr).
\end{split}
\end{equation*}
Therefore, with miraculous cancellations,
\begin{equation*}
\begin{split}
T_2 + T_3 + T_4 + T_5 =& -\frac{2x\cos{(h\log{x})}}{\pi(4+h^2)}
\int_{1}^{\infty} \frac{\sin{\frac{Ty}{x}}}{y^2} dy +
\frac{T}{\pi} \frac{2B\cos{(h\log{x})}}{4+h^2} \\
&+ \frac{T}{\pi} \frac{4h\sin{(h\log{x})}}{(4+h^2)^2} + O\Bigl(T
\bigl(\frac{T}{x}\bigr)^{1/2-\epsilon}\Bigr) +
O\Bigl(\frac{\tilde{h}T}{\log^M{T}}\Bigr).
\end{split}
\end{equation*}
By Lemma \ref{lemma6.1} and $B = -C_0 - \log{2\pi}$,
\begin{equation*}
\begin{split}
F_h(x,T) =& \frac{T}{\pi} \Bigl[\frac{2 \cos{(h\log{x})}}{4 + h^2}
\log{x}\Bigr] - \frac{2T\cos{(h\log{x})}}{\pi(4+h^2)}
\Bigl[\frac{\sin{(T/x)}}{T/x} - ci\bigl(\frac{T}{x}\bigr)\Bigr]
\\
&+ \frac{T}{\pi} \frac{2B\cos{(h\log{x})}}{4+h^2} + O\Bigl(T
\bigl(\frac{T}{x}\bigr)^{1/2-\epsilon}\Bigr) +
O\Bigl(\frac{\tilde{h} T} {\log^{M-2}{T}}\Bigr) \\
=& \frac{T}{\pi} \Bigl[\frac{2 \cos{(h\log{x})}}{4 + h^2}
\log{x}\Bigr] - \frac{2T\cos{(h\log{x})}}{\pi(4+h^2)} \Bigl[1 -
C_0 - \log{\frac{T}{x}} \\
&+ C_0 + \log{2\pi}\Bigr] + O\Bigl(T
\bigl(\frac{T}{x}\bigr)^{1/2-\epsilon}\Bigr) +
O\Bigl(\frac{\tilde{h} T} {\log^{M-2}{T}}\Bigr) \\
=& \frac{T}{2\pi} \log{\frac{T}{2\pi e}} \Bigl[\frac{4
\cos{(h\log{x})}}{4 + h^2}\Bigr] + O\Bigl(T
\bigl(\frac{T}{x}\bigr)^{1/2-\epsilon}\Bigr) +
O\Bigl(\frac{\tilde{h} T} {\log^{M-2}{T}}\Bigr).
\end{split}
\end{equation*}
\section{Sketch for Conjecture \ref{conj2}}

Fix $\alpha > 0$. Let $r(u)$ be an even function which is almost
the characteristic function of the interval $[-\alpha,\alpha]$
with $\hat{r}(\alpha) \ll \frac{1}{\alpha^2}$ (see page 87 of
[\ref{C1}] for detail construction). We use Conjecture \ref{conj1}
to compute the right hand side of (\ref{0}).
\begin{equation*}
\begin{split}
I =& \int_{-\infty}^{\infty} F_h(\alpha) \hat{r}(\alpha) d\alpha =
2 \int_{0}^{\infty} F_h(\alpha) \hat{r}(\alpha) d\alpha \\
=& 2(1+o(1)) \log{T} \int_{0}^{1} T^{-2\alpha} \hat{r}(\alpha)
d\alpha + 2 \frac{4}{4+h^2} \int_{0}^{1} \alpha \cos{(h \log{T}
\alpha)}
\hat{r}(\alpha) d\alpha \\
&+ 2 \frac{4}{4+h^2} \int_{1}^{\infty} \cos{(h\log{T} \alpha)}
\hat{r}(\alpha) d\alpha + O\Bigl(\frac{1}{A}\Bigr) + o(1) \\
=& \frac{4}{4+h^2} \int_{-\infty}^{\infty} \cos{(h\log{T} \alpha)}
\hat{r}(\alpha) d\alpha - \frac{4}{4+h^2} \int_{-1}^{1} (1 -
|\alpha|) \cos{(h\log{T} \alpha)} \hat{r}(\alpha) d\alpha \\
&+ (1+o(1)) \log{T} \int_{-\infty}^{\infty} T^{-2|\alpha|}
\hat{r}(\alpha) d\alpha + O\Bigl(\frac{1}{A}\Bigr) + o(1) \\
=& \frac{4}{4+h^2} \int_{-\infty}^{\infty} \hat{r}_1(\alpha)
d\alpha - \frac{4}{4+h^2} \int_{-1}^{1} (1 - |\alpha|)
\hat{r}_1(\alpha) d\alpha \\
&+ (1+o(1)) \log{T} \int_{-\infty}^{\infty} T^{-2|\alpha|}
\hat{r}(\alpha) d\alpha + O\Bigl(\frac{1}{A}\Bigr) + o(1) \\
=& \frac{4}{4+h^2} I_1 - \frac{4}{4+h^2} I_2 + (1+o(1)) I_3 +
O\Bigl(\frac{1}{A}\Bigr) + o(1)
\end{split}
\end{equation*}
where $r_1(u) = r(u + \frac{h\log{T}}{2\pi})$. As
$\int_{-\infty}^{\infty} \hat{r}_1(\alpha) d\alpha = r_1(0)$,
$$I_1 = r_1(0) = r\Bigl(\frac{h\log{T}}{2\pi}\Bigr).$$
By $\int f \hat{g} = \int \hat{f} g$, the transform pair and the
definition of $r(u)$,
$$f(t) = \max{(1-|t|,0)},\, \mbox{ } \, \hat{f}(u) =
\Bigl(\frac{\sin{\pi u}}{\pi u}\Bigr)^2,$$
$$I_2 = \int_{-\infty}^{\infty} r_1(u) \Bigl(\frac{\sin{\pi u}}
{\pi u}\Bigr)^2 du = \int_{-\alpha + h\log{T}/(2\pi)}^{\alpha +
h\log{T}/(2\pi)} \Bigl(\frac{\sin{\pi u}} {\pi u}\Bigr)^2 du +
o(1).$$ Similarly, by the transform pair
$$f(t) = e^{-2a|t|},\, \mbox{ } \, \hat{f}(u) = \frac{4a}{4a^2 +
(2\pi u)^2},$$
$$I_3 = \int_{-\alpha}^{\alpha} \frac{4\log^2{T}}{4\log^2{T} +
(2\pi u)^2} du + o(1) = \int_{-\alpha + h\log{T}/(2\pi)}^{\alpha +
h\log{T}/(2\pi)} 1 du + o(1).$$ Therefore,
\begin{equation}
\label{7.1} I = \frac{4}{4+h^2} r\Bigl(\frac{h\log{T}}{2\pi}\Bigr)
+ \int_{-\alpha + h\log{T}/(2\pi)}^{\alpha + h\log{T}/(2\pi)} 1 -
\frac{4}{4+h^2} \Bigl(\frac{\sin{\pi u}} {\pi u}\Bigr)^2 du +
O\Bigl(\frac{1}{A}\Bigr) + o(1).
\end{equation}
Now, the left hand side of (\ref{0}) is
\begin{equation}
\label{7.2} \frac{4}{4+h^2} r\Bigl(\frac{h\log{T}}{2\pi}\Bigr) +
\Bigl(\frac{T}{2\pi} \log{T}\Bigr)^{-1} \mathop{\sum_{0 < \gamma
\neq \gamma' \leq T}}_{|\gamma - \gamma' - h| \leq 2\pi\alpha /
\log{T}} (1 + o(1)).
\end{equation}
Combining (\ref{7.1}) and (\ref{7.2}), we have Conjecture
\ref{conj2} by making $A$ arbitrarily large. The only shaky point
in the above argument is the error analysis. All of these become
rigorous following page $87-90$ of [\ref{C1}].

Tsz Ho Chan\\
Case Western Reserve University\\
Mathematics Department, Yost Hall 220\\
10900 Euclid Avenue\\
Cleveland, OH 44106-7058\\
USA\\
txc50@po.cwru.edu

\end{document}